\pdfoutput=1 %% arXiv admin added this line

\documentclass{amsart}
\usepackage{graphicx}
\usepackage{amssymb}
\usepackage{epstopdf}
\theoremstyle{plain}
\newtheorem{Lem}{Lemma}[section]

\newtheorem{Thm}[Lem]{Theorem}
\newtheorem{Prop}[Lem]{Proposition}
\newtheorem{Cor}[Lem]{Corollary}
\theoremstyle{definition}

\newenvironment{Pf}{\noindent {\em Proof.}}{\hfill\raisebox{5mm}[6mm]{\framebox[2mm][l]{ }}
\smallskip}
\usepackage{color}
\definecolor{purple}{rgb}{0.6,0,0.7}

\begin{document}

\title{Moving faces to other places: Facet derangements}
\author{{Gary Gordon}
\and{Elizabeth McMahon}}
\address{Lafayette College\\
   Easton, PA  18042\\ {\tt mcmahone@lafayette.edu}}
 \thanks{The authors would both like to thank the Isaac Newton Institute, Cambridge, UK and Lafayette College for support during this project.}


\begin{abstract}
Derangements are a popular topic in combinatorics classes.  We study a generalization to face derangements of the $n$-dimensional hypercube.  These derangements can be classified as odd or even, depending on whether the underlying isometry is direct or indirect, providing a link to abstract algebra.  We emphasize the interplay between the geometry, algebra and combinatorics of these sequences, with lots of pretty pictures.
\end{abstract}

\maketitle

\section{Introduction}  \label{S:Intro}
Suppose you have a die sitting on a table.  Instead of just looking at the number on the top of the die, pay attention to the location of the numbers on all  six sides.  Now  you pick it up and roll it so that it occupies the same place on the table it did before.  How many ways could you have done this so that none of the 6 numbers are in the same place?  

\begin{center}
\includegraphics[width=3in]{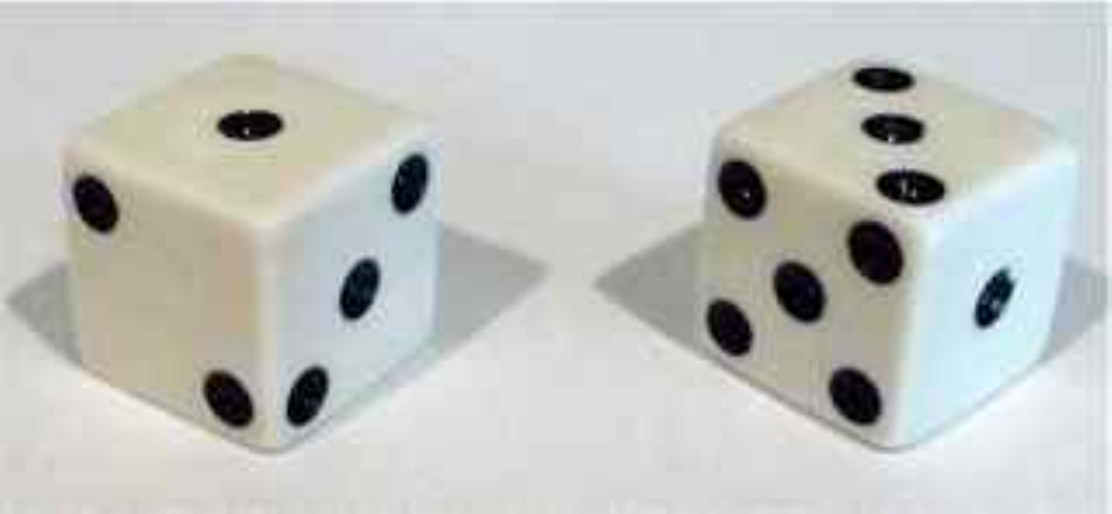}
\end{center}

Rolling a die so that it occupies the same place it did before it was rolled gives a {\it direct isometry} of the cube, i.e., a geometric transformation realizable in 3 dimensions that fixes the cube.  Then our question becomes:

\begin{quote}
How many direct isometries of the cube are {\em derangements} of the faces of the cube?
\end{quote}

How many additional derangements  do you get if you allow yourself to turn your cube inside out (allowing an indirect isometry via a roll through 4-dimensional space)?  Generalizing, if you have an $n$-dimensional hypercube, then the same question makes sense, although it's a bit harder to buy one and actually roll it.

We'll need a few important facts about isometries in $n$ dimensions.  Every isometry can be written as a composition of reflections through hyperplanes.  An isometry is direct if it can be written as a composition of an even number of reflections, otherwise, it is indirect (as with products of transpositions in symmetric groups, this is well-defined).  Further, a direct isometry is orientation preserving, while an indirect isometry reverses orientation. The composition of two reflections is a rotation, provided the hyperplanes corresponding to those reflections intersect.  This will always be the case with the isometries we consider here since we will be considering isometries of regular solids whose centers will always be fixed.

The total number of isometries  of an $n$-cube is  $2^n n!$ -- see Chapters 5 and 7 of Coxeter's classic \cite{cox} for a guide to understanding the geometry of these isometries.  Half of these are direct and half are indirect.  

We now have several questions at hand:
\begin{itemize}
\item  How many of the $2^n n!$ isometries are derangements of the ($n-1$)-dimensional faces (also called facets) of the $n$-cube?  
\item How many of these facet derangements are direct and how many are indirect?  
\item What if our die is not a cube (or hypercube)?  In particular, what are the counts for a die shaped as an $n$-dimensional simplex?
\end{itemize}
The answers will lead us to integer sequences that have been studied before in several contexts, as well as to  two new sequences.    In particular, the sequences associated with the direct and indirect facet-derangements of the $n$-cube do not appear in Sloane's Online Encyclopedia of Integer Sequences (OEIS, \cite{OEIS}).

These problems lie at the intersection of three fields: combinatorics, geometry and algebra.  Our philosophy is mathematically inclusive here.   Derangements are typically studied in combinatorics classes, where they provide  a good example of inclusion-exclusion.  Students need to take a class in abstract algebra to see even and odd permutations.  Finally,   isometries of Euclidean space 
might appear in a geometry class emphasizing transformations.  But no one class typically presents all three topics coherently.

Derangements of faces of different dimensions of the cube have been studied before.  In \cite{cs}, Chen and Stanley derive explicit generating functions that count derangements of the vertices, edges and 2-dimensional faces for an $n$-cube; the number of isometries that {\it do} fix at least one vertex is given by an especially attractive formula:  $(2n-1)!!=(2n-1)(2n-3)\cdots3\cdot 1$.  Chen and Zhang  study the {\it excedance} of a signed permutation   in \cite{chen} and \cite{cz}, and  Chen, Tang and Zhao \cite{ctz}   generalize this approach to derangement {\it polynomials}.  In a different direction,  Shareshian and Wachs give an interesting connection to combinatorial topology in \cite{sw}, where Theorem 6.2 shows that the number of facet derangements of the $n$-cube is the dimension of the reduced homology of the order complex of a certain poset.

Most of the material we present here is known, but our approach seeks to unify the  combinatorics, algebra and geometry.  This paper is organized as follows:  In Section \ref{S:DSimp}, we begin with a treatment of ordinary derangements from combinatorial, geometric and algebraic viewpoints.  
 Section \ref{S:Dcube} generalizes the entire approach from Section \ref{S:DSimp} to the hypercube.  This gives us a `new' combinatorics problem (we call it the {\it coatcheck problem}, generalizing the hatcheck problem for ordinary derangements), and  (what we believe to be) a new recursion for the facet derangements for the $n$-cube.   We also present several known formulas for the number of facet derangements, with proofs that use ideas from combinatorics (inclusion-exclusion) and linear algebra (signed permutation matrices).
 
Section \ref{S:orient} studies the partition of the facet derangements of the ($n-1$)-simplex and the $n$-cube into direct and indirect isometries.  This is where the connections to geometry are deepest, and where some very pretty relations are developed.

The geometry of 3 dimensions gives us a chance to test our geometric intuition, so Section \ref{S:3D} describes these isometries and derangements in some detail.  We conclude by offering a few suggestions for further study in Section \ref{S:future}.

%%%%%%%%%%%%%%%%%%%%%%%%%%%%%%%%%%%

\section{Derangements and   simplices}\label{S:DSimp}
\subsection{Counting ordinary derangements}  One of the standard problems in a combinatorics class asks the following question:
\begin{quote}
{\bf Hatcheck problem:} $n$ people check their hats at the beginning of a party.  How many ways can the hats be returned later   so that no one gets their own hat back?
\end{quote}

Slightly more modern versions might involve returning cell phones to students after an exam, or designing a cryptoquote so that no letter stands for itself (usually called a {\it substitution cypher}) in the coded message.

Permutations of $\{1,2, \dots , n\}$ with no fixed points are called {\it derangements}.  The following formula counting derangments can be found in every combinatorics book in the section that introduces {\it inclusion-exclusion} as a counting technique.  If $d_n$ denotes the number of derangements on an $n$-element set, then 
\begin{equation}\label{Eq:derange}
d_n=n!\sum_{k=0}^n\frac{(-1)^k }{k!}
\end{equation}

From the formula,  we get $d_0 = 1$.  You can give any explanation you like for what happens when no people  check their hats, as long as you get the answer 1.

This formula has a very attractive probabilistic consequence:  If the hats are randomly returned to the party-goers, then the probability that no one receives their own hat is approximately $e^{-1}\sim 36.79\%$.  When $n$ is large, this probability is (essentially) independent of $n$, which is surprising (unless you already know this, in which case it isn't).  The number of derangements for $n \leq 7$ is given in Table \ref{T:derange}.  
Note that the ratio $d_7/7!$ agrees with $e^{-1} \approx 0.367879\dots$ to 4 decimals.

\begin{table}[htdp]
\renewcommand{\arraystretch}{1.2}
\caption{Number of derangements $d_n$ for $n \leq 7.$}
\begin{center}
\begin{tabular}{|c|c|c|c|c|c|c|c|c|} \hline
$n$ & 0 & 1&2&3&4&5&6&7 \\ \hline
$d_n$ & 1 & 0&1&2&9&44&265&1854 \\ \hline
$\displaystyle{\frac{d_n}{n!}}$ \rule[-.35cm]{0cm}{.9cm} & 1 &0 & 0.5& $0.\bar{3}$& 0.375 & $0.3\bar{6}$ &$0.3680\bar{5}$ & $0.367857\dots$ \\ \hline
\end{tabular}
\end{center}
\label{T:derange}
\end{table}%

We end this preliminary discussion with one more formula -- a very useful recursion for $d_n$ that we'll refer to later:

\begin{equation}\label{Eq:derangerecursion}
d_n=(n-1)(d_{n-1}+d_{n-2})
\end{equation}

The proof of this recursion follows from partitioning the derangements into those in which 1 participates in a 2-cycle and those in which 1 is in a longer cycle.  A proof can be found in \cite{br}, for instance. 

\medskip

\subsection{ Derangements and the $(n-1)$-simplex}    Our first  goal is to interpret derangements geometrically, as derangements of the facets (or vertices) of the $(n-1)$-simplex.  Recall the $(n-1)$-simplex is formed by joining $n$ affinely independent points in $\mathbb{R}^{n-1}$, so the 2-simplex is a triangle in the plane, the 3-simplex is a tetrahedron in 3-space, and so on.
 \begin{quote}
{\bf Geometric Interpretation:} A derangement on $\{1, 2, \dots, n\}$ corresponds to an isometry in $\mathbb{R}^{n-1}$ of the regular $(n-1)$-simplex  in which every  one of the $n$ facets is moved.
\end{quote}

Why does this work?  First, the full symmetry group (including both direct and indirect isometries) of the regular $(n-1)$-simplex is the symmetric  group $S_n$.  The group is usually thought of as acting on the $n$ vertices of the simplex, but since every facet has a unique vertex  it does not touch,  a permutation moves a vertex  if and only if it moves the `opposing' facet.  Thus, for our purposes, we will think of $S_n$ as acting on the facets of the $(n-1)$-simplex .

\begin{figure}[h]
\begin{center}
\includegraphics[width=3in]{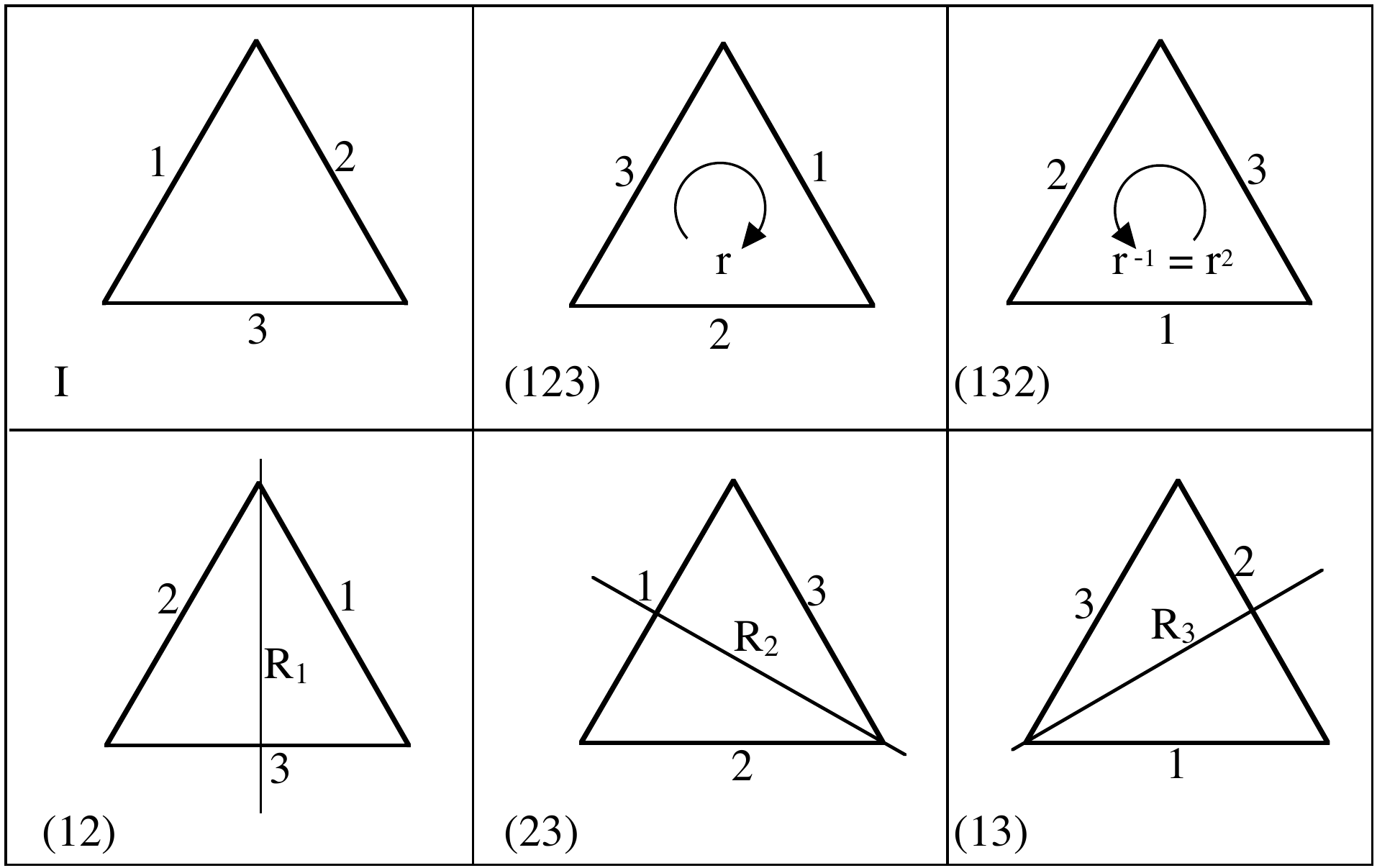}
\caption{Symmetry group $=S_3 \cong D_3$}
\label{F:triangles}
\end{center}
\end{figure}

In 2 dimensions, a simplex is an equilateral triangle in the plane.  There are only two kinds of isometries possible:  rotations (direct) and reflections (indirect).   This group of isometries  is easy to visualize via the six edge permutations   of the triangle.  The correspondence between the isometries and the permutations is given in Fig. \ref{F:triangles}.  

Which of these are derangements?  The two rotations $r=(123)$ and $r^2=(132)$ derange the edges, and the remaining four permutations fix at least one edge, so $d_3=2$, as we saw earlier.

In 3 dimensions, we need to add one more isometry to our tool kit: {\it rotary reflections}.  These can be hard to visualize - they correspond to compositions of reflections and rotations.  
For some fun (and a good exercise in geometric visualization), try to figure out the geometry of the  facet derangements in 3 dimensions.  You should get the following:
\begin{itemize}
\item For the regular tetrahedron, there are 9 derangements: 3 of these are rotations and the remaining 6 are  rotary reflections.
\end{itemize}
Much more information about the geometry of isometries in 3 dimensions -- and some very attractive pictures -- appears in Section \ref{S:3D}.  Feel free to skip ahead.

%%%%%%%%%%%%%%%%%%%%%%%%%%%%%%%%%%%%%%%%

\section {Deranging the facets of a hypercube} \label{S:Dcube}
\subsection{A problem with coats and cubes}
In Section \ref{S:DSimp}, we saw that ordinary derangements correspond to isometries of a regular $(n-1)$-simplex in which no ($n-1$)-dimensional facet is fixed. To generalize this to facet derangements of hypercubes, we need to modify the hatcheck problem:
\begin{quote}
{\bf Couples coatcheck problem} This time, $n$ {\it couples} each check their two coats at the beginning of a party; the attendant puts a couple's coats on a single hanger.  The coats are returned at the end of the party in the following way: when a couple arrives to get their coats, the (lazy) attendant picks an arbitrary hanger and then hands one of those coats to one person in the couple and the other coat to the other (again, arbitrarily).  How many ways can the coats be returned   so that no one gets their own coat back?
\end{quote}

If no one receives their own coat, we'll say we have a {\it c-derangement} (where `c' stands for `coat' or `couple' or `cube'). Evidently, there are two ways the people in a couple could fail to receive their own coats:  Either the pair of coats belonging to that couple was given to another couple, or the couple did receive their own coats, but the coats were swapped between the two partners.

As before, we begin by interpreting this combinatorial problem geometrically.  

 \begin{quote}
{\bf Geometric Interpretation:} A c-derangement  corresponds to an isometry in $\mathbb{R}^n$ of the hypercube in which no ($n-1)$-dimensional facet is fixed.
\end{quote}

Why is this true?  First, note that there are $2^n n!$ ways to return the coats with no restrictions:
\begin{itemize}
\item First, permute the $n$ hangers in $n!$ ways;
\item Next,  hand the coats back to the two members of each couple in $2^n$ ways. 
\end{itemize}

But this is precisely the number of isometries of the $n$-dimensional hypercube:
\begin{itemize}
\item First, permute the $n$ facets around a given vertex in $n!$ ways;
\item Next, swap or don't swap each pair of opposite facets in $2^n$ ways.
\end{itemize}

\begin{figure}[h]
\begin{center}
\includegraphics[width=3in]{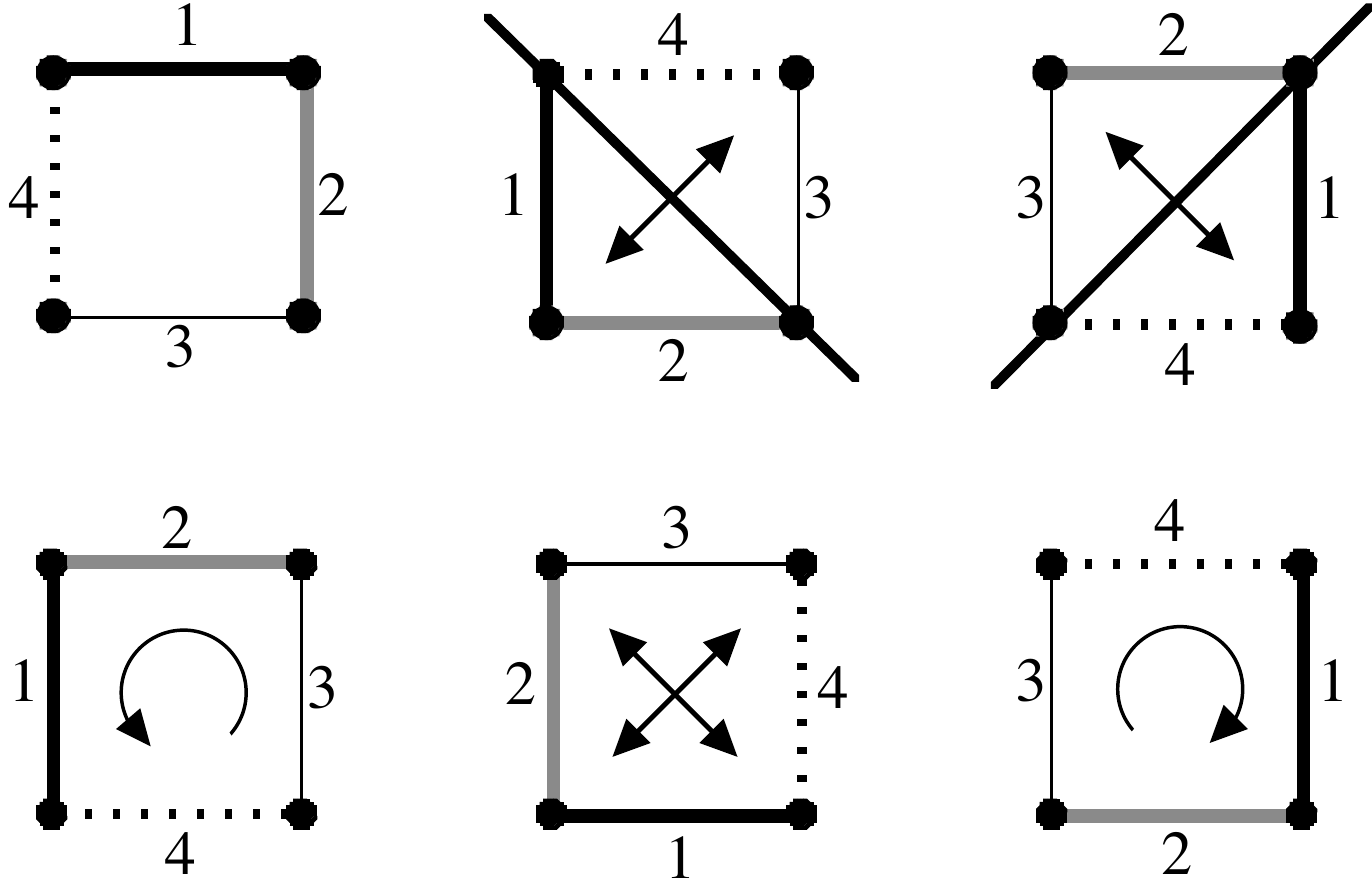}
\caption{The 5 side-derangements of the square.}
\label{F:Sq}
\end{center}
\end{figure}

 Let's introduce some notation:  Let $Q_n$ denote the regular $n$-dimensional hypercube.  The $2n$ people are represented by the $2n$ symbols $1, 1^*, 2, 2^*, \dots, n, n^*$, where $\{i,i^*\}$ is the $i^{th}$ couple.   Then the $2n$ facets of $Q_n$ correspond precisely to the $2n$ people in the coatcheck problem, with the couple $\{i,i^*\}$ corresponding to that pair of opposite facets in the hypercube.  It should now be clear that moving all the facets of $Q_n$ is equivalent to deranging the coats.  (Unlike the connection between  ordinary derangements and the simplex, deranging vertices is {\it not} the same as deranging facets.  In Figure \ref{F:Sq}, you can locate a c-derangement of the square that is not a vertex derangement, and then find a vertex-derangement that's not a c-derangement.)
 
 A 2-dimensional cube is a square, and its facets are the  4  bounding edges.    Among the 8 symmetries of a square (which form the dihedral symmetry group $D_4$), there are 5 c-derangements, pictured in Figure \ref{F:Sq}.  The original square is in the upper left.  The other two squares in the top row arise from reflections through the diagonals of the square, and the 3 squares in the bottom row of the figure arise from  rotations.

\bigskip

As with the simplex, the geometry is fun to explore in 3 dimensions.  We return to this (in great detail) in Section \ref{S:3D}.  For now, though, see if you can determine which isometries of the cube are derangements.  As a check, you should get the following:
\begin{itemize}
\item For the cube, there are 29 derangements: 14 rotations and 15 rotary reflections.
\end{itemize}

\subsection{Some formulas for c-derangements}  We are ready to derive some formulas for the number of c-derangements.  Let $\hat{d}_n$ be the number of c-derangements of $\{1,1^*,\dots,n,n^*\}$.  As in the ordinary derangement case, we can find a formula for $\hat{d}_n$ using inclusion-exclusion.  This formula is listed for sequence A000354 in the OEIS \cite{OEIS}, but the interpretation in terms of c-derangements is not mentioned there.

\begin{Thm} \label{T:NumDer} The number $\hat{d}_n$ of c-derangements of the facets of $Q_n$ is given by 
\begin{equation}\label{Eq:cderange}
\hat{d}_n =  \sum_{k=0}^n (-1)^k {n \choose k} 2^{n-k} (n-k)!
\end{equation}

\end{Thm}

\begin{Pf}  First,  place the $n$-dimensional hypercube so that one vertex is at the origin and the opposite vertex is at the point $(1,1,\dots,1)$.  Then $Q_n$ has $n$ facets with one vertex at the origin, and we need only consider what happens to those facets to determine what happens to all the facets of $Q_n$.  

We use inclusion-exclusion.  The number of c-derangements will be the total number of isometries of $Q_n$, minus the number that fix at least one facet, plus the number that fix at least 2 facets, etc.  

The number of isometries of $Q_n$ that fix at least 0 facets is just the total number of isometries of $Q_n$, that is, $2^n n! $.  To fix $k$ (or more) facets, we choose those $k$ in ${n \choose k}$ ways, and then consider the number of automorphisms of the rest, which is $2^{n-k} (n-k)!$.  This gives the formula.

\end{Pf}

Note that a consequence of this formula is that $\hat{d}_0 = 1$.  Again, you can figure out what that means yourself.
Also notice that this formula guarantees that $\hat{d}_n$ is always odd: every term in Equation \ref{Eq:cderange} is even except for the term corresponding to $k = n$.  We'll see two more proofs of this fact later -- pay attention.

As with ordinary derangements,  we can use this result to get a nice probabilistic interpretation.  Rewriting the above formula as \begin{equation}\label{Eq:cderange2}
\hat{d}_n=2^n n!   \sum_{k=0}^n (-1/2)^k \frac{1}{k!} 
\end{equation}
gives us the following:
\medskip
\begin{quote}
In the coatcheck problem, the probability that no one receives their own coat approaches $e^{-1/2}$ as the number of couples increases.
\end{quote} 

\medskip

Table \ref{T:cderange} gives the number of c-derangements for $n \leq 7$.  When $n=6$, the approximation to $e^{-1/2}\approx 0.606531\dots$ is accurate to 5 decimals.  This series converges faster than the series for $e^{-1}$ that gave the number of (ordinary) derangements, as you would expect from the Taylor series error estimate.
\begin{table}[htdp]
\caption{Number of c-derangements $\hat{d_n}$ for $n \leq 7.$}
\begin{center}
\begin{tabular}{|c|c|c|c|c|} \hline
$n$ &0&1&2&3\\ \hline
$\hat{d_n}$ \rule[-.25cm]{0cm}{.6cm} &1& 1&5&29 \\ \hline
$\displaystyle{\frac{\hat{d_n}}{2^nn!}}$ \rule[-.35cm]{0cm}{.9cm} &1 & 0.5& 0.625 & $0.6041\dots$ \\ \hline \hline
$n$&4&5&6&7 \\ \hline 
$\hat{d_n}$ \rule[-.25cm]{0cm}{.6cm} &233&2329&27,949&391,285 \\ \hline
$\displaystyle{\frac{\hat{d_n}}{2^nn!}}$ \rule[-.35cm]{0cm}{.9cm} & $0.606770\dots$& $0.606510\dots$ &$0.606532\dots$ & $0.606530\dots$ \\ \hline
\end{tabular}
\end{center}
\label{T:cderange}
\end{table}%

There is another formulation for $\hat{d_n}$ that we like.   We can express $\hat{d_n}$ with a sum that uses the ordinary derangements $d_k$.  The proof uses signed permutation matrices, which will be  important for us in Section \ref{S:orient} when we consider the orientation of a c-derangement, i.e., whether the c-derangement corresponds to a direct or an indirect isometry of the hypercube.

\begin{Prop} \label{P:NumDer} The number of c-derangements of the facets of $Q_n$ is \begin{equation}\label{Eq:cderange3}
\hat{d_n} =  \sum_{k=0}^n {n \choose k} 2^k d_k .
\end{equation}

\end{Prop}

\begin{Pf}  Situate $Q_n$ so that all its vertices are located at the points  $(\pm 1,\pm1, \dots,\pm1)$.  Then the centers of the facets are $\pm e_i$, where $e_i$ is the $i^{th}$ standard basis vector in $\mathbb{R}^n$.  Now any isometry of $Q_n$ sends the facet corresponding to $e_i$ to another facet $\pm e_j$.  Thus, an isometry can be identified with a {\it signed permutation matrix}, a permutation matrix with with entries $\pm 1$. (There are clearly $2^n\cdot n!$ such matrices.) 

Now let $A$ be a signed permutation matrix that corresponds to a facet derangement.  Then there are no 1's on the main diagonal of $A$.  (An entry of $-1$ along the main diagonal corresponds to sending a facet to its opposite.)  It remains to count the number of such matrices.

Suppose $A$ has $k$ 0's along its main diagonal.  Then the remaining $n-k$ diagonal entries must be $-1$'s.  We can view the action of $A$ on $Q_n$ as follows:  $A$ fixes the $n-k$ facet pairs  corresponding to the $-1$'s on the diagonal (swapping the facets within each such pair), and $A$ deranges the remaining $k$ facet pairs.  We can choose the $k$ facet pairs to derange in $n \choose k$ ways, then derange those pairs in $d_k$ ways.  But there are two choices for each such pair: If the pair $\{e_i,-e_i\}$ is sent to $\{e_j,-e_j\}$, then either  $e_i \mapsto e_j$ (and so $-e_i \mapsto -e_j$), or  $e_i \mapsto -e_j$ (so $-e_i \mapsto e_j$).

Thus, the total number of c-derangement signed permutation matrices  having $k$ 0's on their diagonals is just  ${n \choose k}2^k d_k$.   The formula follows by summing over $k$.

\end{Pf}

\begin{figure}[h]
\begin{center}
\includegraphics[width=1in]{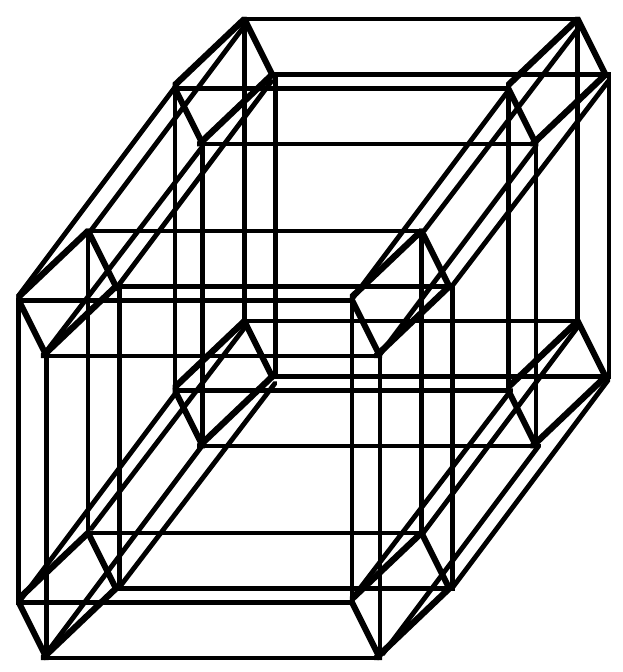}
\caption{A 5-cube with two sets of opposite facets that will be interchanged; the remaining 3 pairs will be deranged.}
\label{F:HF}
\end{center}
\end{figure}

One way to visualize the derivation of  Equation (\ref{Eq:cderange3}) is  suggested by Figure \ref{F:HF}.  In the figure, assume that two pairs of opposite facets of a 5-dimensional cube are interchanged, so $i \leftrightarrow i^*$ and $j \leftrightarrow j^*$, corresponding to two $-1$'s on the main diagonal of the signed permutation matrix.  Then the two facet pairs $(i,i^*)$ and $(j,j^*)$ are `shrunk,' leaving a 3-cube.  That cube will then have its pairs of opposite faces deranged in $d_3$ ways (since we don't allow opposite faces to be interchanged), and then pairs of opposite sides can be oriented in $2^3$ ways.

A generalization of Equation (\ref{Eq:cderange3}) appears in \cite{ss}, where  Spivey and Steil define the {\it rising $k$-binomial transform} of a sequence $a_n$ to be 
\[ r_n = \left\{ \begin{array}{c l}
\sum_{i=0}^n {n \choose i} k^i a_i  &   \mbox{if } k \neq 0   \\
a_0  &   \mbox{if } k = 0.   \\
\end{array}  \right.\] 
Then the sequence $\hat{d}_n$ is the rising 2-binomial transform of the  sequence $d_k$.  The explicit connection to c-derangements does not appear in that paper, however.

The following recursion generalizes the recursion for ordinary derangements given in Equation (\ref{Eq:derangerecursion}).  The recursion appears to be new.   It also gives an inductive proof that $\hat{d_n}$ is always odd (our second proof of this fact).

\begin{Prop}\label{P:crecursion} For $n \geq 2, \hat{d_n}$ satisfies:
\begin{equation}\label{Eq:crecursion} 
\hat{d_n}=(2n-1)\hat{d}_{n-1}+(2n-2)\hat{d}_{n-2}
\end{equation}

\end{Prop}
\begin{Pf}
This proof is modeled on the proof of Equation~(\ref{Eq:derangerecursion}).  That proof partitions derangements into 2 classes:  those in which $1$ is interchanged with some $k>1$ and those in which $1$ is in a cycle of length greater than 2.  For c-deragements, we have an additional case arising from the fact that a pair $(i,i^*)$ can be fixed but swapped and additional complications because there are multiple ways to send $(i,i^*)$ to $(j,j^*)$.

There are three cases to consider for c-derangements.

\begin{description}
\item[Case 1] The pair $(1,1^*)$ is fixed by the c-derangement.  In this case, $1$ and $1^*$ must be swapped and the remaining $n-1$ pairs are c-deranged.  This gives a total of $\hat{d}_{n-1}$ c-derangements.  (This case does not appear in the proof of the ordinary derangement recursion (\ref{Eq:derangerecursion}).)
\item[Case 2] The pair $(1,1^*)$ is swapped with some other pair $(k,k^*)$. There are $n-1$ choices for the $(k,k^*)$ pair, and then there are 4 ways to swap the pairs $(1,1^*)$ and $(k,k^*)$:
(1) $1 \leftrightarrow k $, $ 1^* \leftrightarrow  k^* $; or 
(2) $1 \leftrightarrow k^*$,  $1^* \leftrightarrow  k $; or 
(3) $1 \rightarrow k \rightarrow 1^* \rightarrow k^* \rightarrow 1$; or (4) $1 \rightarrow k^* \rightarrow 1^* \rightarrow k \rightarrow 1$.

The remaining $n-2$ pairs can be c-deranged in $\hat{d}_{n-2}$ ways.  This gives a total of $4(n-1)\hat{d}_{n-2}$ c-derangements.
\item[Case 3] 1 is in a cycle of length greater than 2.  If $(1,1^*) \rightarrow (k,k^*)$, then by assumption, $k \neq 1$ and $(k,k^*) \not \mapsto (1,1^*)$.  There are $n-1$ choices for $k$, and there are 2 ways to send $(1,1^*)$ to $(k,k^*)$.  Since $(k,k^*) \mapsto (i,i^*)$, where $(i,i^*) \neq (1,1^*)$, we can simply replace $(k,k^*)$ by either $(1,1^*)$ or $(1^*,1)$, depending on how $(1,1^*)$ mapped to $(k,k^*)$ in the original c-derangement.  This gives us a c-derangement on $n-1$ pairs with the additional restriction that the pair $(1,1^*)$ is not mapped to itself.  (This restriction is a consequence of the fact that we can't send $(k,k^*)$ to $(1,1^*)$.)

How many of the $\hat{d}_{n-1}$ derangements require that $(1,1^*)$ not be sent to itself?  The number of c-derangements of $n-1$ pairs in which one pair {\em is} sent to itself is  $\hat{d}_{n-2}$ (this is just Case 1 for the remaining $n-1$ pairs).  Thus, the final total in this case is $2(n-1)(\hat{d}_{n-1} - \hat{d}_{n-2})$ derangements.

\end{description}

These three cases give that $\hat{d}_n = \hat{d}_{n-1} + 4(n-1)\hat{d}_{n-2} + 2(n-1)(\hat{d}_{n-1} - \hat{d}_{n-2}) = (2n-1)\hat{d}_{n-1}+(2n-2)\hat{d}_{n-2}$, as desired.

\end{Pf}

Notice that $\hat{d_2}$ and the values you get from  the recursion depend on the (somewhat disturbing) fact that $\hat{d_0}=1$. 

%%%%%%%%%%%%%%%%%%%%%%%%%%%%%%%%%%%%%%%%%%%%%

\section{The parity of derangements and c-derangements}\label{S:orient}

\subsection {Counting direct and indirect (ordinary) derangements}
As we've seen, isometries can be either direct or  indirect, depending on whether they can be expressed as a product of an even or odd number of reflections, resp.  The direct isometries of a simplex correspond to even permutations of $\{1,2,\dots,n\}$, and the indirect ones correspond to odd permutations.

For the facet derangements of a simplex, here's what we have seen so far:  For the  triangle,  both derangements are rotations, and so are direct isometries.  For the regular tetrahedron, there are 9 derangements:   3 rotations (direct) and 6 rotary reflections (indirect).  See Section \ref{S:3D} for a full explanation for the tetrahedron.

We let $e_n$  denote the number of direct derangements of the facets of an $(n-1)$-simplex and $o_n$ denote the number of indirect derangements.   Table \ref{Ta:evenodd} gives the values of $d_n,e_n$ and $o_n$ for $n \leq 7$.

\begin{table}[htdp]
\caption{Number of even and odd derangements  for $n \leq 7.$}
\begin{center}
\begin{tabular}{|c|c|c|c|c|c|c|c|} \hline
$n$ & 1&2&3&4&5&6&7 \\ \hline
$d_n$& 0&1&2&9&44&265&1854 \\ \hline
$e_n$& 0&0&2&3&24&130&930 \\ \hline
$o_n$& 0&1&0&6&20&135&924 \\ \hline \hline
$e_n-o_n$& 0&$-1$&2&$-3$&4&$-5$&6 \\ \hline
\end{tabular}
\end{center}
\label{Ta:evenodd}
\end{table}%

A very casual glance at the last row of the table should suggest a very pretty theorem.  This pattern is not new; for example,  you can find the sequences $e_n$ and $o_n$ in  the OEIS \cite{OEIS}.   A bijective proof using a conjugation argument can be found in \cite{ch}, and the recent paper  \cite{mr} offers two more proofs based on an analysis of excedances of permutations.  The proof we give is inductive and more in line with what we have done to this point.

\begin{Thm}\label{T:evenodd}  
Let $e_n$ and $o_n$ be the respective number of even and odd derangements of $\{1,2,\dots,n\}$.  Then $e_n-o_n=(-1)^{n-1}(n-1)$.
\end{Thm}
\begin{Pf}  We give a recursive procedure for computing $e_n$ and $o_n$.  The proof then follows by induction.  The key observation is the familiar recursion given in Equation (\ref{Eq:derangerecursion}): $d_n=(n-1)(d_{n-1}+d_{n-2})$.

Adapting this recursion to include the parity of the derangements is easy:  If 1 is in a transposition, we just remove that transposition.  If 1 is not in a transposition, then we just remove 1 from the cycle it participates in.  In either case, this process changes even permutations to odd ones, and vice versa.

We let $e_n'$ denote the even derangements in which the transposition $(1 r)$ appears for some $r>1$ and $e_n''$ denote the other even derangements, and define $o_n'$ and $o_n''$ similarly.  Then $e_n=e_n'+e_n'' $ and $ o_n=o_n'+o_n'' $.  This gives the following recursive relations:
$$e_n'=(n-1)o_{n-2}  \hskip.2in e_n''=(n-1)o_{n-1} \hskip.2in o_n'=(n-1)e_{n-2}  \hskip.2in o_n''=(n-1)e_{n-1}$$

Then
\begin{eqnarray*}
e_n-o_n&=&(e_n'+e_n'')-(o_n'+o_n'') \\
&=&(n-1)((o_{n-1}+o_{n-2})-(e_{n-1}+e_{n-2})) \\
&=& (n-1)((o_{n-1}-e_{n-1})+(o_{n-2}-e_{n-2})) \\
& = & (n-1)((-1)^{n-1}(n-2) + (-1)^{n-2}(n-3)) \\
& = & (-1)^{n-1}(n-1).
\end{eqnarray*}

The formula now follows by induction, using the initial values from Table \ref{Ta:evenodd}.

\end{Pf}

\begin{Cor} Let $e_n$ and $o_n$ be the number of even and odd derangements, resp.  Then $\displaystyle{e_n=\frac{d_n+(-1)^{n-1}(n-1)}{2}}$ and $\displaystyle{o_n=\frac{d_n+(-1)^{n}(n-1)}{2}}$.
\end{Cor}

We don't believe Theorem \ref{T:evenodd} is as well-known as it should be.  As we remarked earlier, we believe this has to do with the fact that students study derangements in combinatorics classes, but they study permutation groups in algebra classes.

\subsection{Counting direct and indirect c-derangements}
How can we determine whether an isometry of $Q_n$ is direct or indirect?  We know that direct isometries are the product of an even number of reflections, but in this setting, it can be confusing to interpret reflections using permutation notation.  For instance, consider the following c-derangement of the cube in 3-dimensions:

\begin{itemize}
\item First, reflect in a plane that contains two edges of the cube to get the  facet permutation $(1, 2)(1^*,2^*)$.
\item Then reflect in a plane parallel to facets 3 and $3^*$:  $(3, 3^*)$.
\end{itemize}

The result is a c-derangement whose facet cycle structure is $(1,2)(1^*,2^*)(3,3^*)$.  This {\it looks} like a product of an odd number of transpositions, but this isometry is {\it direct} (and so, in this case, it's a rotation).  The problem with our algebraic representations of the isometries is the redundancy in representing the reflection $(1,2)(1^*,2^*)$.

We can work around this problem by finding a more concise way to represent an isometry.  For our example, the representation of the first isometry has the following redundancy:  If facets 1 and 2 are swapped, then the isometry also {\it must} swap $1^*$ and $2^*$.  Thus, we can recode the isometry as follows:  $(1,2)(3,3^*).$  Similarly, we can recode $(1,2^*)(1^*,2)$ as $(1,2^*)$ (of course, we could also recode it as $(1^*,2)$ if we wanted).

For readers who know more group theory, this recoding can be made precise by using the group of isometries of the hypercube, called the {\it hyperoctahedral} group (it's the same as the isometry group of the dual solid - the {\it hyperoctahedron} or the $n$-dimensional {\it cross-polytope} \cite{cox}).  This group is isomorphic to $\mathbb{Z}_2^n \rtimes S_n$, the semidirect product of the normal abelian 2-group $\mathbb{Z}_2^n$ and the symmetric group $S_n$, where the action of $S_n$ on $\mathbb{Z}_2^n$ is conjugation.

From a geometric perspective, we can think of this group as acting on the facets of the hypercube as follows:\begin{itemize}
\item Situate the hypercube so that the origin $\bf{0}$ is one of its vertices   and the facets incident to $\bf{0}$  are labeled $1, 2, \dots, n$.
\item Use elements in $S_n$ to permute the facets surrounding $\bf{0}$.
\item Finally, use reflections in the normal subgroup $\mathbb{Z}_2^n$ to move $\bf{0}$ to some other vertex $v$.
\end{itemize}

We saw this procedure earlier in the discussion at the beginning of Section \ref{S:Dcube} that explained the link between the coat problem and c-derangements.  Any permutation in $S_n$ can be written as a product of transpositions, and you can show that interchanging two facets adjacent to $\bf{0}$ can be realized as a single reflection. %%
Then the first factor $\mathbb{Z}_2^n$ is generated by $n$ orthogonal reflections normal to the $n$ coordinate axes, and $S_n$ is generated by reflections corresponding to transpositions.    More information about this action can be found in \cite{h}, for example.

We can now use this to get the `right' way to represent an isometry of the cube as a product of transpositions.   Our goal is to write any isometry of the cube as a product of terms like $(i,i^*)$ followed by terms of the form $(i,j)$.  (We use the convention that multiplication is done right to left.)   For example, the rotation  $(12)(1^*2^*)(33^*)$ would be expressed more simply as $(33^*)(12)$ (since if 1 and 2 are switched, then $1^*$ and $2^*$ must also be switched).  More complicated but still doable, $(1,2^*)(1^*,2) = (1,1^*)(2,2^*)(1,2)(1^*,2^*)$ which becomes $(1,1^*)(2,2^*)(1,2)$ in our shorthand notation.  

Using this representation for the group, we see that an isometry is direct if and only if it is a product of an even number of transpositions (written in the shorthand way), which is what we were after.

An even better way to deal with the direct-indirect issue here is to use the signed permutation matrices  we used in the proof of Prop. \ref{P:NumDer}.  Then we can tell whether an isometry is direct or indirect by evaluating its determinant:
\begin{quote}
If $A$ is a signed permutation matrix, then $\det(A)=1$ precisely when $A$ corresponds to a direct isometry, and $\det(A)=-1$ if $A$ corresponds to an indirect isometry.
\end{quote}

\begin{table}[htdp]
\begin{center}
\begin{tabular}{|c|c|c|} \hline
Reflection & Transposition & Matrix operation \\ \hline
Normal to coordinate axis & $(k,k^*)$ & Multiply row $k$ by $-1$ \\ \hline
Not normal to coordinate axis & $(ab)(a^*b^*)$ & Swap rows $a$ and $b$ \\ \hline
Not normal to coordinate axis  & $(ab^*)(a^*b)$ & Swap rows $a$ and $b$ \\
 & & and multiply each by $-1$  \\ \hline
\end{tabular}
\end{center}
\label{T:smatrix}
\end{table}%

In Section \ref{S:Dcube} we gave the number of  direct and indirect c-derangements for the square and cube.  Here's what we asserted:
\begin{itemize}
\item Square:  $\hat{d_2}=5$, with 3 direct and 2 indirect c-derangements.
\item Cube:  $\hat{d_3}=29$, with 14 direct and 15 indirect c-derangements.
\end{itemize}

The next theorem shows that the difference between the number of direct and indirect c-derangements is always $\pm 1$.  This will follow from the construction of  a bijection between the odd and even c-derangements {\it except for  central inversion}.   Central inversion  is the isometry that sends every point $(x_1,\dots,x_n)$ to its antipodal point $(-x_1,\dots,-x_n)$.  For the square, you can see the central inversion in Figure \ref{F:Sq} in the center of the bottom row --   a 180$^\circ$ rotation indeed sends every point to its antipode.  Thus, in dimension 2, central inversion is a direct isometry.  In Section \ref{S:3D}, we will see that central inversion in dimension 3 is indirect (chemists call this an {\it improper reflection}).  This pattern continues:  the matrix corresponding to  central inversion is $-I_{n\times n}$, i.e., the diagonal matrix with all entries $-1$.  Its determinant is $\pm 1$, depending on the parity of $n$, the dimension of our hypercube.  Thus, central inversion (which is always a c-derangement) is direct in even dimensions and indirect in odd dimensions.

\begin{Thm}\label{T:evenoddcube}
Let $\hat{e}_n$ and $\hat{o}_n$ denote the number of direct and indirect  c-derangements of $Q_n$  resp.  Then $\hat{e}_n-\hat{o}_n=(-1)^n.$
\end{Thm}
\begin{Pf}   To keep track of which c-derangements are direct and which are indirect, we use signed permutation matrices.  We will construct a bijection between the direct and indirect c-derangments that are not central inversion.  So, assume $A$ is a signed permutation matrix corresponding to a c-derangement that is not central inversion. $A$ has at least 2 non-zero entries off the main diagonal, so we can find $i$ so that  $a_{i,i}=0$ but $a_{m,m}\neq 0$ for all $m<i$.  Let  $j$ be the column corresponding to the non-zero entry in row $i$, so   $j$ is the unique index with  $a_{i,j}\neq 0$.

We know that $a_{i,j}=1$ or  $-1$.  For one of those choices, the determinant will be 1 and for the other, the determinant will be $-1$ (since $a_{i,j}$ is the only non-zero entry in that row, and multiplying a row by $-1$ multiplies the determinant by $-1$).  Since both choices are c-derangements, this gives us the bijection between the direct and indirect c-derangements (where we have removed central inversion from consideration).  Putting central inversion back in  gives us the formula, since central inversion is direct if and only if $n$ is even.

\end{Pf}

Theorem \ref{T:evenoddcube} gives a third proof\footnote{Three proofs mean it's really really really true.} that $\hat{d_n}$ is always odd: $\hat{d_n}=2\hat{e}_n \pm 1$.
In Table \ref{Ta:evenoddcube}, we list the first 7 values of the sequences $\hat{e}_n$ and $\hat{o}_n$.

\begin{table}[htdp]
\caption{Number of even and odd c-derangements  for $n \leq 7.$}
\begin{center}
\begin{tabular}{|c|c|c|c|c|c|c|c|} \hline
$n$ & 1&2&3&4&5&6&7 \\ \hline
$\hat{d}_n$& 1&5&29&233&2329&27,949&391,285 \\ \hline
$\hat{e}_n$& 0&3&14&117&1164& 13,975& 195,642 \\ \hline
$\hat{o}_n$& 1&2&15&116&1165& 13,974& 195,643  \\ \hline
\end{tabular}
\end{center}
\label{Ta:evenoddcube}
\end{table}%

\begin{Cor}\label{C:evenoddcderange}  Let $\hat{e}_n$ and $\hat{o}_n$ denote the number of direct and indirect  c-derangements of $Q_n$  resp.  Then $\displaystyle{ \hat{e}_n=\frac{\hat{d_n}+(-1)^n}{2}}$ and $\displaystyle{ \hat{o}_n=\frac{\hat{d_n}+(-1)^{n+1}}{2}}$.

\end{Cor}

The  sequences $\hat{e}_n$ and $\hat{o}_n$ appear to be `new' in the sense they are not (currently) in the OEIS \cite{OEIS}.

%%%%%%%%%%%%%%%%%%%%%%%%%%%%%%%%%%%%%%

\section{Geometry of 3-dimensional derangements and c-derangements}\label{S:3D}  Time for some fun.  If, like most humans,  you have trouble visualizing objects in 4 or more dimensions, then you can start by looking at some low dimensional examples, and then try to generalize.  This sounds easier than it is, but, for the simplex and hypercube symmetries, it's extremely valuable.  An interesting discussion about the power of analogy appears in Section 7-1 of \cite{cox}.

\subsection{Derangements of the tetrahedron}
For a bounded solid in $\mathbb{R}^3$, there are three kinds of isometries we need to consider:  rotations, reflections and {\it rotary reflections}.  Rotary reflections are operations that correspond to combining a rotation and a reflection.  It is possible to have a rotary reflection where neither the rotation nor the reflection involved is itself a symmetry of the object; when this happens, the rotary reflection is called {\it irreducible}.

\begin{figure}[htbp]
\begin{center}
\includegraphics[width=4.3in]{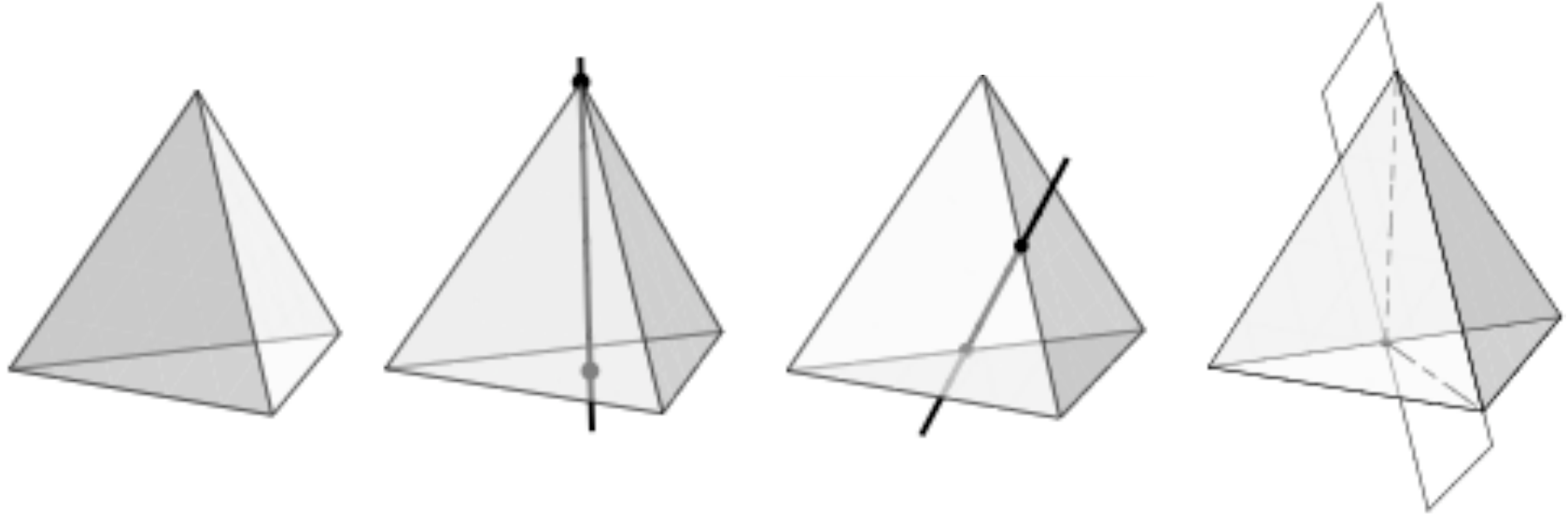}
\newline (a) Identity \hskip.3in (b) $120^{\circ}$ Rotation  \hskip.2in (c) $180^{\circ}$ rotation \hskip.2in (d) Reflection
\caption{A tetrahedron, a $120^{\circ}$ rotation,  a $180^{\circ}$ rotation and a reflection.}
\label{F:tetra2}
\end{center}
\end{figure}

The tetrahedron has 24 symmetries corresponding to the 4! permutations of $S_4$.  As isometries, these can be broken into the following categories. 

\begin{figure}[htbp]
\begin{center}
\includegraphics[width=4in]{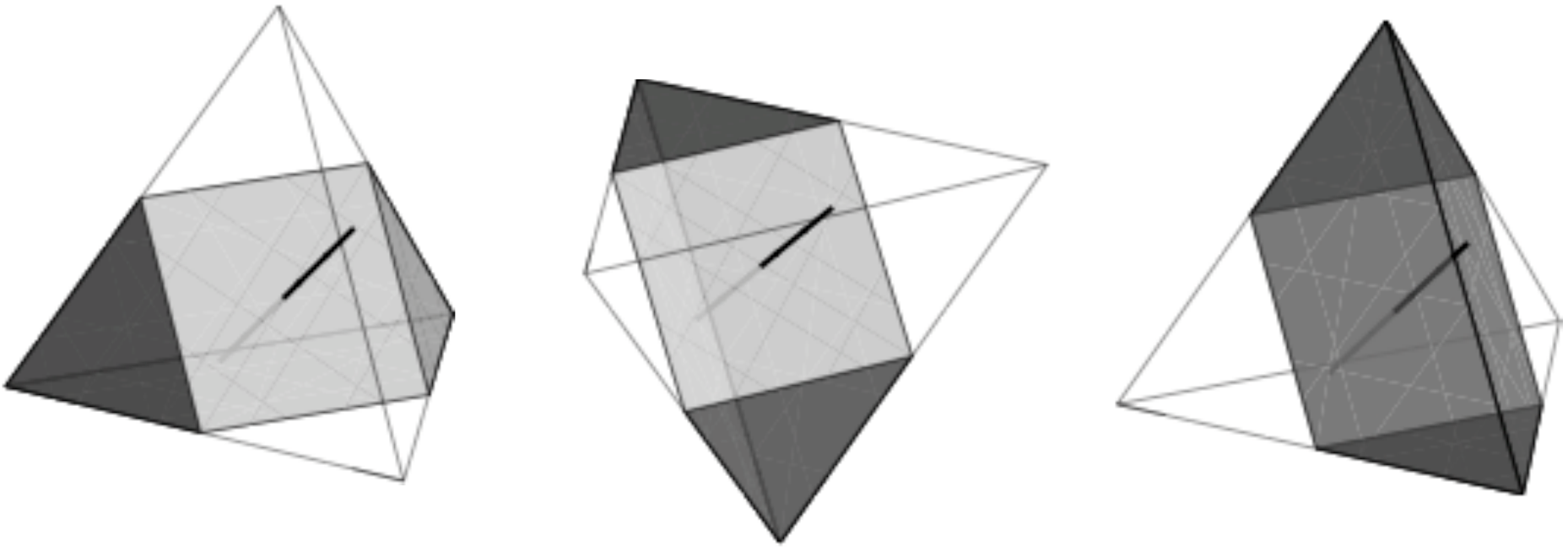}
\newline \hskip.15in Original \hskip.75in  $90^{\circ}$ rotation \hskip.65in Reflection
\caption{A rotary reflection  in two steps:  First rotate (the tetrahedron does not match up), then reflect.}
\label{F:tetra1}
\end{center}
\end{figure}

\renewcommand{\labelitemii}{$\circ$}
\begin{itemize}
\item Direct isometries
\begin{itemize}
\item The identity I, which is clearly not a derangement;
\item The 8 rotations with axis through a  vertex and the center of the opposite facet as in Figure \ref{F:tetra2}(b). These order 3 rotations are through angles of $120^{\circ}$ or $240^{\circ}$ and correspond to face permutations of the form  $(abc)$.    These are not derangements;
\item The 3  $180^{\circ}$  rotations with axis through the centers of two opposite edges as in Figure \ref{F:tetra2}(c).  They have face permutations  $(ab)(cd)$ and so are derangements.
\end{itemize}

\item Indirect isometries
\begin{itemize}
\item The 6 reflections with mirror plane containing one edge of the tetrahedron as in Figure \ref{F:tetra2}(d).  These have face permutations of the form $(ab)$, so they aren't derangements;
\item The 6 rotary reflections, with face permutation form $(abcd)$ as in Figure \ref{F:tetra1}.   These are derangements.
\end{itemize}
\end{itemize}
 
It's worth investigating more closely how the rotary reflections operate.   The following step-by-step guide should help you create any of the 6 rotary reflections on your own tetrahedron.
\begin{enumerate}
\item Cut the tetrahedron through the plane containing the centers of four edges (take all  edges except  two that do not share a vertex).  This cuts the tetrahedron into two congruent pieces that get glued back along a  square.  (You now have a tricky little  puzzle.) 
\item Rotate the entire tetrahedron $90^{\circ}$ about an axis perpendicular to the square cross-section.  This is {\em not} an isometry of the tetrahedron.
\item Reflect the tetrahedron through the plane containing the square.  This is also not an isometry of the tetrahedron.
\end{enumerate}

You didn't really need to cut your tetrahedron into two pieces after all, but you do need to identify the square cross-section.  The rotation involved is `half' of one of our $180^{\circ}$ edge-rotations from Figure \ref{F:tetra2}(c). The resulting rotary reflection is irreducible.  As mentioned above, this rotary reflection induces the face permutation $(abcd)$, a 4-cycle.  This shouldn't be too surprising, as the square from the cross-section meets each of the four faces, and we are rotating around this square.

\subsection{c-derangements of the cube}  
\begin{figure}[h]
\begin{center}
$\begin{array}{c@{\hspace{.6in}}c}
\includegraphics[width=2in]{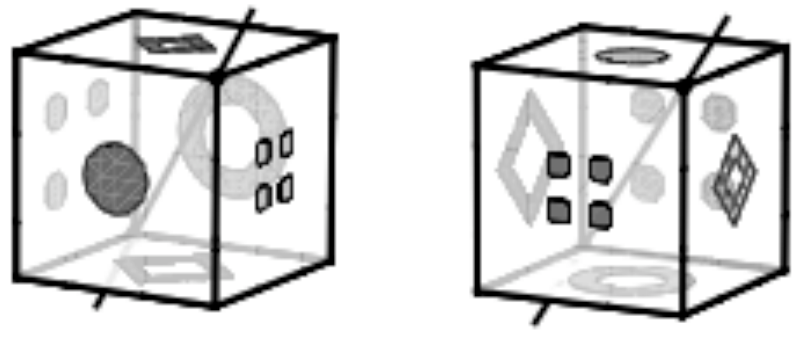} &
\includegraphics[width=2.15in]{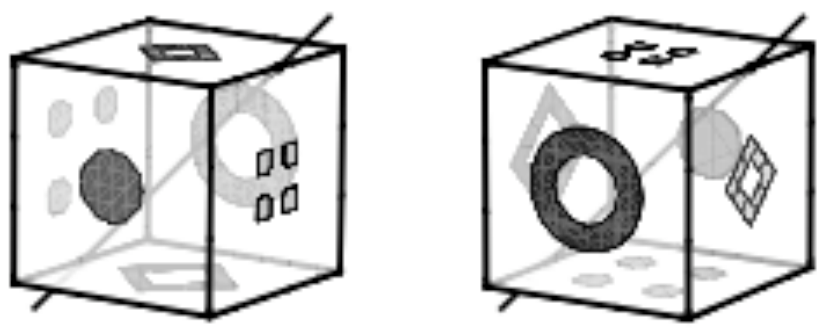} 
\end{array}$
\hspace{1.2in}(a) Rotate through opposite vertices  \hspace{.2in} (b) Rotate through opposite edge centers 
\caption{The direct face-derangements of the cube.}
\label{F:C1}
\end{center}
\end{figure}

\begin{figure}[h]
\begin{center}
$\begin{array}{c@{\hspace{.2in}}c@{\hspace{.2in}}c}
\includegraphics[width=1.48in]{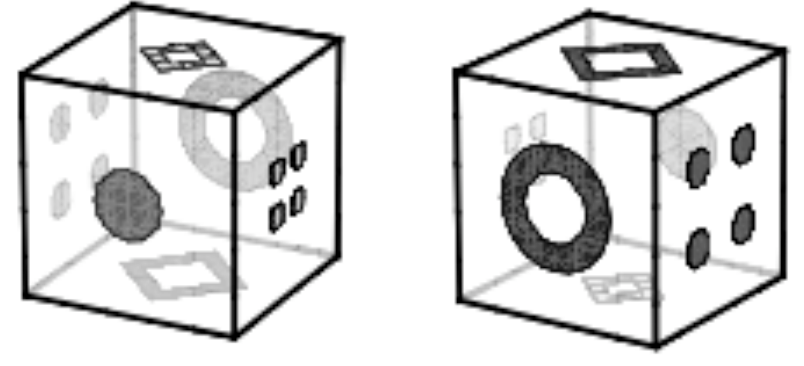} &
\includegraphics[width=1.5in]{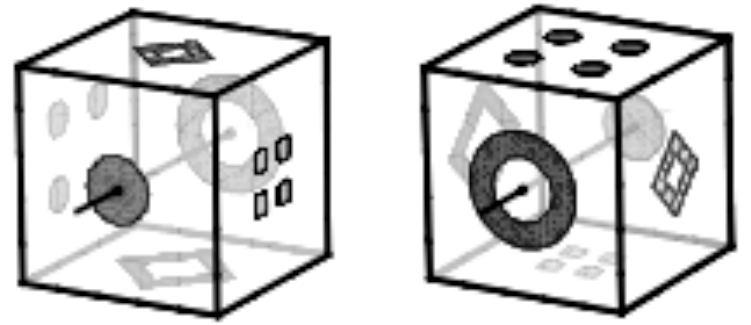} &
\includegraphics[width=1.6in]{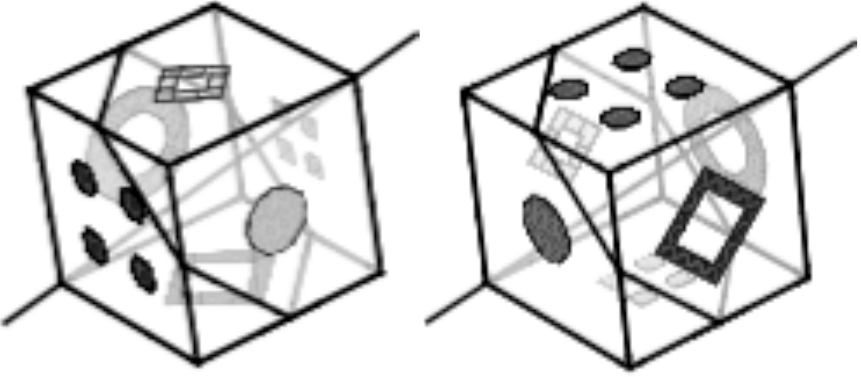}
\end{array}$
\hspace{1.5in}(a) Central inversion  \hspace{.1in} (b) Axis through side centers
 \hspace{.1in} (c) Diagonal axis
\caption{The 15 indirect face-derangements of the cube.}
\label{F:C2}
\end{center}
\end{figure}

Now  we turn our attention to the ordinary 3-dimensional cube.  
There are $2^33! = 48$ isometries of a cube:
\begin{itemize}
\item Direct isometries
\begin{itemize}
\item The identity, which is not a c-derangement;
\item Eight  rotations of 120$^\circ$ and 240$^\circ$ through the 4 pairs of opposite vertices as pictured in \ref{F:C1}(a).  These are c-derangements;
\item  Six $180^{\circ}$ rotations through the centers of opposite edges as pictured in \ref{F:C1}(b).  These are also c-derangements;
\item Nine   rotations through the centers of opposite faces.  These are not c-derangements.
\end{itemize}
\item Indirect isometries
\begin{itemize}
\item Nine reflections: three through the plane containing the centers of four parallel edges and six through  opposite edges.  None of these are c-derangements;
\item Fifteen rotary reflections: central inversion (pictured in Figure \ref{F:C2}(a)); six reducible ones through an axis through centers of opposite sides  (pictured in Figure \ref{F:C2}(b)) and eight irreducible ones through an axis through two opposite vertices  (pictured in Figure \ref{F:C2}(c)).  These are c-derangements.
\end{itemize}
\end{itemize}

Note that 14 of the 29 c-derangements are direct, while 15 are indirect.  In both Figure \ref{F:C1} and Figure \ref{F:C2},  opposite facets have the same pattern, but  in different sizes.  Notice also that the direct isometries preserve the orientation of the cube (so you perform these isometries with a die), while the indirect ones reverse the orientation of the cube (so you'd need a mirror-image die to perform them).

The indirect c-derangements are all rotary reflections, which  are, again, a bit harder to picture.  Central inversion always corresponds to a c-derangement (see  Figure \ref{F:C2}(a)), but it's an indirect isometry in 3-dimensions.  It can be realized either as the composition of 3 mutually perpendicular reflections or by a rotary reflection in a variety of ways, as we will see shortly.  

We now examine the 14 remaining rotary reflections in a bit more detail.   In Figure \ref{F:C2}(b), you can see the rotary reflection formed as a composition of a rotation of 90$^\circ$ or 270$^\circ$ around an axis through the center of two opposite sides followed by a reflection through the plane perpendicular to that axis at its midpoint.  There are 3 of those axes so there are 6 such rotary reflections, all of which are reducible. (Notice that if you do this same rotary reflection with a rotation of 180$^\circ$, you get central inversion).  
 
\begin{figure}[h]
\begin{center}
\includegraphics[width=4in]{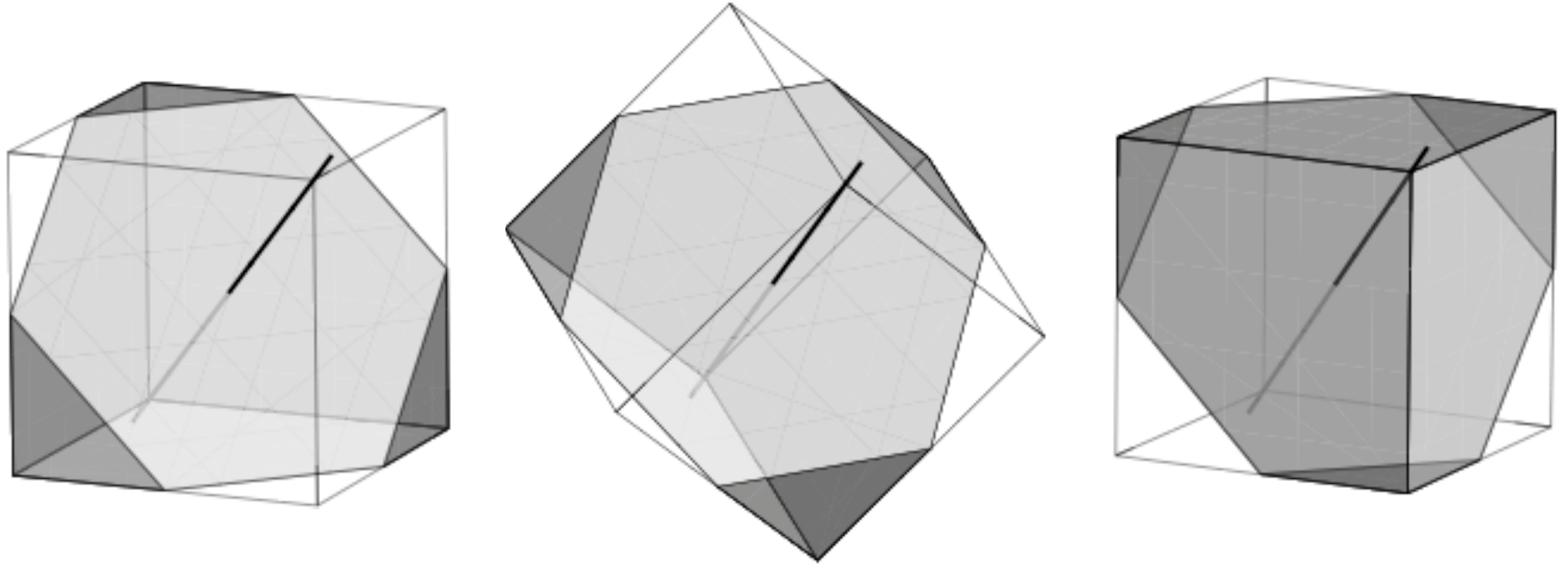}
\caption{An irreducible rotary reflection of the cube.}
\label{F:CRR}
\end{center}
\end{figure}

 In Figure \ref{F:C2}(c), you can see the last kind of rotary reflection, which is irreducible: a rotation of 60$^\circ$ or 300$^\circ$ around an axis through two opposite vertices followed by a reflection through the plane perpendicular to that axis. Since there are 4 pairs of opposite vertices, this gives 8 distinct isometries. (Once again, central inversion arises from a rotation of 180$^\circ$ followed by reflection.)  These can be quite challenging to picture, since the initial rotation does not align the cube with itself.  Figure \ref{F:CRR} breaks an irreducible rotary reflection down to help with visualization. 

 There are two more interesting aspects of this kind of rotary reflection.  First, it gives a 6-cycle on the faces of the cube.  Second, the  plane perpendicular to the axis of rotation cuts the cube into two congruent pieces, forming (as a cross section) a regular hexagon that passes through each of the 6 sides of the cube.  It is well worth locating these hexagons on an actual cube!  Showing students the two pieces (after having them guess about what kinds of polygons could be formed as cross sections when slicing a cube) usually elicits some surprised looks.

%%%%%%%%%%%%%%%%%%%%%%%%%%%%%%%%%%%%%%

\section {Suggestions for future study}\label{S:future}
We conclude with a few ideas for projects that can help solidify some of the ideas from this paper.
\begin{enumerate}
\item In the proof of Theorem \ref{T:evenoddcube}, signed permutation matrices give a bijection between the direct and indirect  isometries (excluding central inversion).  As a warm-up, find this bijection explicitly for the 14 direct and 14 indirect (again, excluding central inversion)
 c-derangements of the cube.  (To do this, you will need to number the six faces of your cube $1,1^*,2,2^*,3,3^*$ and then use this numbering to refer to your c-derangements.)
 
\item Study the vertex, edge and face derangements for the remaining Platonic solids.  We've covered the faces (and the vertices) for the tetrahedron,  and we've done the faces of the  cube (and therefore, the vertices of its dual, the octahedron).  We haven't considered the faces of the octahedron, and we've completely ignored  the icosahedron and its dual, the dodecahedron.  Determining the vertex, edge or facet  derangement numbers   for these solids  is a good exercise in geometric visualization. 

The symmetry group of the icosahedron is $A_5 \times \mathbb{Z}_2$, so there are 120 isometries to consider.  Counting the  direct and indirect vertex, edge and face derangements is also a good exercise.  For extra credit, describe the 45 rotary reflections explicitly in this case.
\item There are 6 regular solids in 4 dimensions:  The 4-simplex, the hypercube   the  hyperoctahedron,  the   24-cell, the  120-cell, and   the 600-cell.  All of the same questions make sense here:
\begin{itemize}
\item Find the number of vertex, edge, 2-dimensional and 3-dimensional face derangement numbers for the 24-cell and the 120-cell.  (By duality, derangements for the 120-cell and the 600-cell are the same, where the vertex derangements for the 120-cell are  the same as the 3-dimesional facet derangements for the 600-cell, and so on.)
\item For each class of derangements, count the  direct and indirect isometries, as above.
\end{itemize}

The 24-cell has 1152 isometries and the 120-cell (and 600-cell) has 14,400.  In dimensions 5 and higher, there are only 3 regular solids:  The $n$-simples, the $n$-cube its dual, the $n$-dimensional hyperoctahedron.  See \cite{h} for descriptions of these symmetry groups.

\item Extend this approach to general root systems.  In particular, the exceptional root systems $E_6, E_7$ and $E_8$ are very important families of vectors that display lots of symmetry.  How many isometries move all the vectors in the root system?  Much more about root systems can be found in \cite{h}.
\end{enumerate}

%****************************************

{\em Acknowledgments.}  The authors wish to thank Ethan Berkove and Derek Smith for helpful comments on the manuscript.

%**************************************


\begin{thebibliography}{99} 


\bibitem{br} R. Brualdi, {\it Introductory Combinatorics},  Pearson Prentice Hall, New Jersey (2004).

\bibitem{ch} R. Chapman, {\it An involution on derangements},  Disc. Math {\bf 231} (2001), 121-122.


\bibitem{chen} W. Chen, {\it Induced cycle structures of the hyperoctahedral group},  SIAM J. Disc. Math {\bf 6} (1993), 353-362.


\bibitem{cs} W. Chen and R. Stanley, {\it Derangements of the $n$-cube},  Disc. Math {\bf 115} (1993), 65-75.


\bibitem{ctz}
W. Chen, R. Tang, and A.  Zhao, {\it Derangement Polynomials and Excedances of Type $B$},  Elec. J. Comb. {\bf 16(2)} R15 (2009).

\bibitem{cz} W. Chen and J. Zhang, {\it The skew and relative derangements of type $B$,} Elec. J. Comb. {\bf 14} N24 (2007).


\bibitem{cox} H.S.M. Coxeter, {\it Regular Polytopes}, Dover, New York (1973).

\bibitem{h}
J. Humphreys, {\it Reflection groups and Coxeter groups}, Cambridge Univ. Press, Cambridge, (1990).
\bibitem{mr} R. Mataci and F. Rakotandrajao, {\it Exceedingly deranging!},  Adv. Appl.. Math {\bf 30} (2003), 177-188.

\bibitem{r} J. B. Remmel, {\it A note on a recursion for the number of derangements},  Europ. J. Combinatorics {\bf 4} (1983), 371-374.

\bibitem{sw}
J. Shareshian and M. Wachs, {\it Poset homology of Rees products, and q-Eulerian polynomials,} arXiv:0812.0779 (2008).

\bibitem{OEIS} N. J. A. Sloane, {\em The Online Encyclopedia of  Integer Sequences}, 2007, published electronically at http://www.research.att.com/~njas/sequences/A000354.

\bibitem{ss} M. Spivey and L. Steil, {\it The $k$-Binomial Transforms and the Hankel Transform},  J. Integer Seq., {\bf 9} (2006), Article 06.1.1.

\end{thebibliography}
 \end{document}